\documentclass{amsart}

\newtheorem{theorem}{Theorem}

\author{Melvyn B. Nathanson}
\address{Department of Mathematics, 
Lehman College (CUNY), 
Bronx, NY 10468} 
\email{melvyn.nathanson@lehman.cuny.edu}

\title{ON A DIOPHANTINE EQUATION OF M. J. KARAMA}


\date{\today}

\begin{document}

\maketitle

\begin{abstract} 
For every positive integer $n$, 
the infinite family of positive integral solutions of the diophantine equation  
$x^n - y^n = z^{n+1}$ is constructed.  
\end{abstract}

\section{The equation $x^n - y^n = z^{n+1}$}
In a recent paper, M. J.  Karama~\cite{mune16} 
studied the diophantine equation $x^2 - y^2 = z^3$,
and conjectured that the diophantine equation 
$x^3 - y^3 = z^4$ has no solution in positive integers.  
A standard reference for diophantine equations is the book by Mordell~\cite{mord69}, 
but this very interesting equation is not discussed there.

We shall prove that, for every positive integer $n$, the diophantine equation 
$x^n - y^n = z^{n+1}$ has infinitely many positive integral solutions.

\section{Powerful triples}

The triple $(a,b,c)$ of positive integers is called \emph{$n$-powerful} 
if $a > b$ and  $c^{n+1}$ divides $a^n - b^n$.  
Define the function 
\begin{equation}                                            \label{tabc}
t_n(a,b,c) = \frac{a^n - b^n }{c^{n+1} }.
\end{equation}
The triple $(a,b,c)$ of positive integers is $n$-powerful 
if and only if $t_n(a,b,c)$ is a positive integer.  
The triple $(a,b,c)$ is \emph{relatively prime} if $\gcd(a,b,c) =1$, 
where $\gcd$ is the greatest common divisor.

\begin{theorem}       \label{DE:theorem:2}
Let $n$ be a positive integer.
If $(a,b,c)$ is an $n$-powerful triple with $t = t_n (a,b,c)$, 
then the triple of positive integers 
\begin{equation}                \label{xyz}
(x, y, z) =  \left( at, bt, ct \right)
\end{equation}
is a solution of the diophantine equation
\begin{equation}         \label{DE-n}
x^n - y^n = z^{n+1}.
\end{equation}
Moreover, there is a one-to-one correspondence between positive integral 
solutions of~\eqref{DE-n} 
and relatively prime $n$-powerful triples.  
\end{theorem}

For example, if $a$ and $b$ are positive integers with $a > b$,  then 
the triple $(a,b,1)$ is $n$-powerful with  $t = t_n(a,b,1) = a^n - b^n$, 
and so 
\begin{equation}   \label{ab1}
(x, y, z) =  \left(  at ,  b t, t \right)  =  \left( a (a^n - b^n),  b  (a^n - b^n),  a^n - b^n \right)  
\end{equation}
is a positive integral solution of~\eqref{DE-n}.  
Moreover,  
\[
\left( a (a^n - b^n),  b  (a^n - b^n),  a^n - b^n \right) 
=  \left( a_1 (a_1^n - b_1^n),  b_1  (a_1^n - b_1^n),  a_1^n - b_1^n \right)
\] 
if and only if $a=a_1$ and $b = b_1$.  It follows that,   
for every $n$, the diophantine equation~\eqref{DE-n} 
has infinitely many solutions.  

Different $n$-powerful triples $(a,b,c)$ can generate identical 
solutions to~\eqref{DE-n}.  
For example, for every positive integer $n$, the triple $(8,4,2)$ is $n$-powerful 
with $t = 2^{2n-1} - 2^{n-1}$, and produces the solution 
\[
(x,y,z) = \left( 2^{2n+2} - 2^{n+2}, 2^{2n+1} - 2^{n+1}, 2^{2n} - 2^{n} \right) 
\]
of the diophantine equation~\eqref{DE-n}.  
The triple $(4,2,1)$ is also $n$-powerful 
with $t = 2^{2n} - 2^{n}$, and produces exactly 
the same solution of~\eqref{DE-n}.

\begin{proof}
Let  $(a,b,c)$ be an $n$-powerful triple with $t = t_n(a,b,c)$. 
Defining $(x,y,z)$ by~\eqref{xyz}, we obtain 
\begin{align*}
x^n - y^n 
& = (at)^n - (bt)^n \\
& = a^n \left(\frac{a^n - b^n}{c^{n+1}} \right)^n - b^n 
\left(\frac{a^n - b^n}{c^{n+1}} \right)^n \\
& = \left( a^n - b^n \right)  \left(\frac{a^n - b^n}{c^{n+1}} \right)^n \\
& =   \left(\frac{a^n - b^n}{c^n} \right)^{n+1} \\
& =   \left(c  \left(\frac{a^n - b^n}{c^{n+1}} \right) \right)^{n+1} \\
& = (ct)^{n+1}  \\
& = z^{n+1}.
\end{align*}
Thus, $(x,y,z)$ solves~\eqref{DE-n}.  

Let $(a,b,c)$ be an $n$-powerful triple with $t =  t_n(a,b,c)$, 
and let $d$ be a common divisor of $a$, $b$, and $c$.  
The relatively prime triple $(a/d,b/d,c/d)$ is $n$-powerful because 
\begin{align*}
t' & = t_n(a/d,b/d,c/d) = \frac{ (a/d)^n - (b/d)^n}{ (c/d)^{n+1}} \\
& = d \left(\frac{a^n - b^n }{c^{n+1} }\right) = d\  t_n(a,b,c) \\ 
& = dt
\end{align*}
is a positive integer.  
The solution of equation~\eqref{DE-n} constructed from $(a/d,b/d,c/d)$ is
\[
(x,y,z) = ( (a/d)t' ,(b/d)t',(c/d)t' ) = (at,bt,ct)
\]
which is also the solution constructed from $(a,b,c)$.

If $(x,y,z)$ is a positive integral solution of the diophantine equation~\eqref{DE-n}, 
then  $(x,y,z)$ is an $n$-powerful triple with $t_n(x,y,z) = 1$.  
Let $d = \gcd(x,y,z)$, and define $(a,b,c) = (x/d,y/d,z/d)$.  
It follows that $(a,b,c)$ is an $n$-powerful triple with 
$t_n(a,b,c) = dt_n(x,y,z) = d$, and that $(x,y,z)$ is the solution of~\eqref{DE-n} 
produced by $(a,b,c)$. 
Thus, every positive integral solution of~\eqref{DE-n} 
can be constructed from a relatively prime $n$-powerful triple.

Let $(x,y,z)$ be a positive integral solution of~\eqref{DE-n}, 
and let $(a,b,c)$ and $(a_1,b_1,c_1)$ be relatively prime $n$-powerful triples 
that produce  $(x,y,z)$.  We must prove that $(a,b,c) = (a_1,b_1,c_1)$.  

If $t = t_n(a,b,c)$ and $t' = t_n(a_1,b_1,c_1)$, then 
\[
(x,y,z) = (at ,bt, ct) = (a_1 t', b_1t', c_1t').
\]
If $d = \gcd(t,t')$, then $t/d$ and $t'/d$ are  positive integers.    
The equation $x = at =a_1t'$ implies that $a(t/d) =a_1(t'/d)$, 
and so $t/d$ divides $a_1( t'/d)$.   
Because  $t/d$ and $t'/d$ are relatively prime, it follows that $t/d$ divides $a_1$,
and $a_1 = A(t/d)$ for some positive integer $A$.  Therefore, 
\[
a \left( \frac{t}{d} \right) =  a_1\left(\frac{ t'}{d}  \right) 
=   A \left( \frac{t}{d}  \right)\left( \frac{t'}{d}  \right)
\]
and $a = A(t'/d)$.  Similarly, there exist positive integers $B$ and $C$ such that 
$b = B (t'/d)$, $b_1 = B( t/d)$,  $c = C (t'/d)$, and $c_1 = C(t/d)$.  
Because $t'/d$ is a common divisor of $a$, $b$, and $c$, and because $\gcd(a,b,c)=1$,
it follows that $t'/d=1$ and so $a = A$, $b = B$, and $c = C$.  
Because $\gcd(a_1,b_1,c_1) = 1$, we also have $a_1 = A$, $b_1 = B$, and $c_1 = C$.
Therefore, $(a,b,c) = (A,B,C) = (a_1, b_1, c_1)$.  
This completes the proof.  
\end{proof}

\section{Open problems}

A Maple computation  
produces 39 positive integral solutions of $x^3 - y^3 = z^4$ 
with $x \leq 5000$.  There are 35 relatively prime $3$-powerful triples of the 
form $(a,b,1)$, and the following four  relatively prime $3$-powerful 
triples $(a,b,c)$ with $c > 1$:
\begin{align*}
t_3(71,23,14) & = 9 \\
t_3(39,16,7) & = 23 \\
t_3(190,163,21) & = 13 \\
t_3(103, 101, 7) & = 26.
\end{align*}
How often is a difference of cubes  divisible by a nontrivial fourth power?  
More generally, how often is a difference of $n$th powers  divisible 
by a nontrivial $(n+1)$st power?

It would also be interesting to know, for positive integers $n$ and $k \geq 2$,   
the positive integral solutions of the diophantine equation 
\[
x^n - y^n = z^{n+k}. 
\]

The Beal conjecture~\cite{maud97} states that if $k,\ell, m$ are integers 
with $\min(k,\ell,m) > 2$ and if $x,y,z$ are positive 
integers such that 
\[
x^k - y^\ell = z^m
\]
then $\gcd(x,y,z) > 1$.  
Does the Beal conjecture hold for the diophantine equation $x^n - y^n = z^{n+1}$?

\def\cprime{$'$} \def\cprime{$'$}
\providecommand{\bysame}{\leavevmode\hbox to3em{\hrulefill}\thinspace}
\providecommand{\MR}{\relax\ifhmode\unskip\space\fi MR }
\providecommand{\MRhref}[2]{%
  \href{http://www.ams.org/mathscinet-getitem?mr=#1}{#2}
}
\providecommand{\href}[2]{#2}

\begin{table}  
\begin{center}
\caption{Solutions of $x^3 - y^3 = z^4$ for $x \leq 5000$ with the associated 
relatively prime $3$-powerful triples $(a,b,c)$ and $t_3 = t_3(a,b,c)$.  
An asterisk (*) indicates a solution with $c > 1$.}

\begin{tabular}{|r | r| r| r| r| r| r |}  \hline
$x$ & $y$ & $z$ & $a$ & $b$ & $c$ & $t_3$ \\  \hline\hline
14 & 7 & 7 & 2& 1& 1& 7\\ \hline
57 & 38 & 19 &3 & 2&1 & 19\\ \hline
78 & 26 & 26 & 3& 1& 1& 26\\ \hline
148 & 111 & 37 & 4& 3& 1& 37\\ \hline
 224 & 112 &  56 & 4 & 2 & 1 & 56 \\ \hline
252  & 63 & 63 &4 & 1& 1& 63\\ \hline
305  & 244 & 61 & 5& 4&1 &61 \\ \hline
490 & 294 & 98 & 5 & 3 & 1 & 98 \\ \hline 
546  & 455 &  91& 6& 5&1 &91 \\ \hline
585  & 234 & 117 & 5& 2& 1 & 117\\ \hline
620  & 124 & 124 & 5& 1 & 1&124 \\ \hline
*639  & 207 & 126 & 71& 23& 14& 9\\ \hline
889  & 762 & 127 & 7& 6& 1& 127\\ \hline
*897 & 368 & 161 & 39& 16& 7&23 \\ \hline
912  & 608 &152  & 6& 4&1 & 152\\ \hline
1134  &  567 & 189 & 6& 3& 1& 189\\ \hline
1248 & 416 & 208 &6 &2 &1 & 208\\ \hline
1290 & 215 & 215 &6 & 1& 1&215 \\ \hline
1352  & 1183 & 169 &8 &7 & 1& 169\\ \hline
1526  & 1090 & 218 & 7& 5&1 & 218\\ \hline
1953  & 1116 & 279 &7 &4 &1 & 279\\ \hline
1953  &1736  & 217 & 9&8 &1 & 217\\ \hline
2212 & 948 & 316 & 7 & 3 & 1 & 316 \\ \hline
2345 & 670 & 335 &7 &2 & 1& 335\\ \hline
2368 & 1776 & 296 & 8 & 6 & 1 & 296 \\ \hline 
2394 & 342 & 342 & 7 & 1 & 1 & 342 \\ \hline 
*2470 & 2119 & 273 & 190& 163& 21&13 \\ \hline 
*2678 & 2626 & 182 & 103 & 101 & 7& 26\\ \hline 
2710 & 2439 & 271 &10 & 9 & 1 & 271 \\ \hline 
3096 & 1935 & 387 & 8 & 5 & 1 & 387  \\ \hline 
3474 & 2702 & 386 & 9  &  7 & 1 & 386\\ \hline
3584 & 1792 & 448 & 8 & 4 & 1 & 448 \\ \hline
3641 & 3310 & 331 & 11 & 10 & 1 & 331 \\ \hline
3880 & 1455 & 485 & 8 &3  &1 &  485 \\ \hline. 
4032   & 1008  &  504 &  8 & 2  & 1&  504 \\ \hline
4088   & 511 &  511 & 8  & 1  &1 &  511 \\ \hline
4617   & 3078 &  513 &  9 & 6  & 1&  513 \\ \hline
4764   & 4367 & 397  & 12  &  11 & 1&  397 \\ \hline
4880   & 3904 & 488  & 10  &  8 & 1 & 488  \\ \hline
\end{tabular}
\end{center}
\end{table}


\begin{thebibliography}{1}

\bibitem{mune16}
{M. J.  Karama}, \emph{{Using summation notation to solve some diophantine
  equations}}, {Palestine Journal of Mathematics} \textbf{5} (2016), 155--158.

\bibitem{maud97}
R.~D.~Mauldin, \emph{A generalization of Fermat's Last Theorem:  
The Beal Conjecture and Prize Problem},  
Notices Amer. Math. Soc. \textbf{44} (1997), 1436--1437.


\bibitem{mord69}
L.~J. Mordell, \emph{{Diophantine Equations}}, Academic Press, London, 1969.

\end{thebibliography}
\end{document}